\documentclass{siamltex}
\usepackage[english]{babel}
\usepackage{graphicx}
\usepackage{textcomp}
\usepackage{amsmath}
\usepackage{amssymb}

\begin{document}

\title{Nonexistence of smooth solutions for the general compressible Ericksen -- Leslie equations in three dimensions}

\author{Tudor S. Ratiu \thanks{Section de Math\'ematiques and Bernoulli Center, \'Ecole Polytechnique F\'ed\'erale de Lausanne, CH-1015 Lausanne,
Switzerland ({\tt tudor.ratiu@epfl.ch}). Partially supported by Swiss NSF grant 200020-126630
and by the government grant of the Russian Federation for
support of research projects implemented by leading
scientists, Lomonosov Moscow State University under the
agreement No. 11.G34.31.0054.}
\and Olga S. Rozanova \thanks{Department of Mechanics and Mathematics, Moscow
State University,  119992 Moscow, Russia} ({\tt rozanova@mech.math.msu.su}). Corresponding author. Partially supported by RFBR Project Nr. 12-01-00308  and by the government grant of the Russian Federation for
support of research projects implemented by leading
scientists, Lomonosov Moscow State University under the
agreement No. 11.G34.31.0054.}

\maketitle

\begin{abstract}We prove that the smooth solutions to the Cauchy problem
for the compressible general three-dimensional Ericksen--Leslie
system modeling nematic liquid crystal flow with
conserved mass, linear momentum, and dissipating total energy,
generally lose  classical smoothness within a finite time.
\end{abstract}



\section{System, known results, and main problem}
\label{1}

The liquid crystal state is often viewed  as an intermediate state
between liquid and solid. The molecules possess none or partial
positional order but display a preferred  orientation. The
nematic phase is the simplest among all liquid crystal phases and is
closest to the liquid state. The molecules float around as in a
liquid, but have the tendency to align along a preferred direction
due to their orientation. The hydrodynamic theory of liquid
crystals, due to Ericksen and Leslie, was developed in the 1960's
\cite{Ericksen1961, Ericksen1962, Leslie1966, Leslie1968}. The full
Ericksen--Leslie system consists of the following equations (cf.
\cite{{Ericksen1987}, {Leslie1968}, {Leslie1979}, {Lin_Liu2000},
{GBRa_2009}}):
\begin{equation}
\partial_t \rho+{\rm div} (\rho u)=0,\label{(1.1)}
\end{equation}
\begin{equation}
\partial_t(\rho u)+{\rm Div} (\rho u \otimes u)+ \nabla
p=-\partial_j\left(\rho \frac{\partial W}{\partial d_j},\nabla
d\right)+{\rm Div}(\sigma),\label {(1.2)}
\end{equation}
\begin{equation}
J\rho \frac{D}{Dt}\omega=(h-g) \times d,\quad \omega=d\times
\frac{D}{Dt} d, \label{(1.4)}
\end{equation}
\begin{equation}
p=A\rho^\gamma,\quad A>0,\,\gamma>1.\label{(1.5)}
\end{equation}

Equations (\ref{(1.1)}) –- (\ref{(1.5)}) are given in
$\mathbb{R}\times\mathbb{R}^3$ and represent the conservation of mass,
linear momentum, and angular momentum, respectively, with the
anisotropic feature of liquid crystal materials exhibited in
(\ref{(1.4)})  and its nonlinear coupling in (\ref{(1.2)}) (cf.
\cite{Leslie1968, Lin_Liu2000}).

Here, we consider the flow of a compressible isentropic  material: $\rho$ is the fluid density, $p$ is a pressure, $J$
is a positive inertial constant, $\gamma$ is an adiabatic
exponent, $ u = (u_1,u_2, u_3)^T$ is the velocity,
$d = (d_1, d_2, d_3)^T$ is the orientational order parameter
representing the macroscopic average
of the molecular directors, and $\frac{D}{Dt}:=
\frac{\partial}{\partial t}+ (u, \nabla )$ denotes the
material derivative.

The notations $A= \frac12 (\nabla u + (\nabla u)^T)$,
$\Omega = \frac12 (\nabla u - (\nabla u)^T)$,
$N= \frac{D}{Dt}d - \Omega d$ represent the rate of strain
tensor, the skew-symmetric part of the strain rate (the vorticity
of the fluid), the material derivative of $d$
(transport of center of mass), and the rigid rotation part of
the changing rate of the director by the fluid vorticity,
respectively.

The Oseen-Z\"ocher-Frank free energy functional $W$ for the
equilibrium configuration of a unit director field for a nematic
crystal is given by the sum of the splay, the bend, and the twist, i.e,
\begin{equation}
W = K_1 \frac12 ({\rm div}\, d)^2 + K_2  \frac12 |d \times ({\rm
curl}\, d)|^2 + K_3  \frac12 (d\cdot{\rm curl}\, d)^2\ge
0,\label{free_energy}
\end{equation}
and the vector field
\begin{equation}\label{h}
h=\rho \frac{\partial W}{\partial d}-\partial_i \left(\rho
\frac{\partial W}{\partial d,_i} \right)
\end{equation}
is the molecular field.

The kinematic transport  $g$ of the director $d$  is defined by:
\begin{equation}\label{g}
 g_i = \lambda_1 N_i + \lambda_2 d_jA_{ji},
 \end{equation}
  and represents
the effect of the macroscopic flow field on the microscopic
structure. The material coefficients $\lambda_1$ and $\lambda_2$
reflect the molecular shape (Jeffrey's orbit) and how slippery
the particles are in the fluid.

The stress tensor $\sigma$ has the following form:
\begin{equation}\label{sigma}
\sigma_{ij} =
 \mu_1 d_k A_{kp} d_p d_i d_j + \mu_2 N_i d_j + \mu_3 d_i N_j +
\mu_4 A_{ij} + \mu_5 A_{ik} d_k d_j + \mu_6 d_i A_{jk} d_k .
\end{equation}
 The independent coefficients $\mu_1, ...,
\mu_6,$ which may depend on the material and temperature, are called
Leslie coefficients.

The following relations are frequently introduced in the literature
(cf. \cite{{Leslie1968}, {Leslie1979}}):
\begin{equation}
 \lambda_1 = \mu_2 - \mu_3,\quad \lambda_2 = \mu_5
- \mu_6, \label{(1.15)}
\end{equation}
\begin{equation}
 \mu_2 + \mu_3 = \mu_6 - \mu_5.
\label{(1.16)}
\end{equation}
 The identities (\ref{(1.15)})
are necessary conditions to satisfy the equation of motion
identically (cf. \cite{Leslie1968}, Section 6). The identity
(\ref{(1.16)}) is called Parodi's relation (cf. \cite{Parodi1970}),
which is derived from Onsager's reciprocal relations expressing the
equality of certain relations between flows and forces in
thermodynamic systems out of equilibrium (cf. \cite{Onsager1931}).
Under the assumption of Parodi's relation, we see that the dynamics
of an incompressible nematic liquid crystal flow involve five
independent Leslie coefficients.

As was proved in \cite{GBRa_2009}, $|d(0)|=1$ implies
$|d(t)|=1$ and
$(d(0), \omega(0))=0$ implies $(d(t), \omega(t))=0 $ for all time.

When $d$ is a constant vector field, the system
(\ref{(1.1)}) -- (\ref{(1.5)}) becomes the compressible
Navier-Stokes equations describing the motion of a
compressible fluid.

Since the mathematical structure of the Ericksen--Leslie system is
quite complicated, the study of the full model presents several
mathematical difficulties. Almost all existing investigations were restricted
to  simplified versions.
For the incompressible model in \cite{Lin1989}, Lin introduced a
simplification of the general Ericksen--Leslie system that keeps
many of the mathematical difficulties of the original system by
using  a Ginzburg--Landau approximation to relax the nonlinear
constraint $d \in \mathcal{S}^2.$ Namely, instead of the
restriction $|d| = 1$ he added  the penalty term
$\frac{1}{\epsilon^2} (|d|^2 -1)^2$ to the free energy functional
$W$ and  neglected the inertial constant $J$. Later in
\cite{Lin_Liu1995}, Lin and Liu showed the global existence of weak
solutions and smooth solutions for that approximation. In
\cite{Shkoller_2002}, a very simple proof of local well-posedness for
this coupled system was provided using a contraction mapping
argument. It was  proved that this system is globally well-posed and
has compact global attractors in 2D. For more results on the
approximative system see \cite{ {Lin_Liu1996}, {Lin_Liu_Zhang2007},
{Liu_Shen_Yang2007}, {Sun_Liu2009}, {FFRS2012}, {CR_2012}}. Further,
taking the limit of $\epsilon \to  0$ in the penalty term, one can
get
 a coupling between the compressible Navier--Stokes equations and a
transported heat flow of harmonic maps into ${\mathcal S}^2$
\cite{Lin_Liu2000}. It is a macroscopic continuum description of the
evolution for  liquid crystals of nematic type under the influence
of both the flow field $u$, and the macroscopic description of the
microscopic orientation configurations $d$ of (rod-like) liquid
crystals. Recently, Hong \cite{Hong2011} and Lin-Liu-Wang
\cite{Lin_Liu_Wang2010} showed independently the global existence of
weak solution of an incompressible model  in two dimensional space.
Moreover, in \cite{Lin_Liu_Wang2010}, the regularity of solutions,
except for a countable set of singularities whose projection on the
time axis is a finite set, has been obtained (see also
\cite{Hong_Xin2012}). In \cite{Wang2011}, Wang established a global
well-posedness theory for the incompressible liquid crystals for
rough initial data, provided that $\|u_0\|_{BMO^{-1}} +
\|d_0\|_{BMO} < \epsilon_0$ for some $\epsilon_0>0$. In
\cite{DQS_2012},  regularity and uniqueness for solutions to density
dependent nematic liquid crystals systems  in the Ginzburg--Landau
approximation were established for a bounded domain: in 2D the
system has a global classical solution, in contrast to the 3D case.
A family of exact solutions with finite energy to the incompressible
liquid crystals in two dimensions was constructed in
\cite{Dong_Lei_2012}.

Concerning the compressible 3D case, local existence and uniqueness
of strong solutions for the coupling between the compressible
Navier-Stokes equations and a transported heat flow of
harmonic maps was proved (see \cite{Huang_Wang_Wen2012}),
provided that the initial data $\rho_0, u_0, d_0$ are
sufficiently regular and satisfy a natural compatibility
condition. A criterion for possible breakdown of such a local
strong solution at finite time was given in terms of blow up
of the $L^\infty$-norms of $\rho$ and $\nabla
d$. Alternative blow-up criteria were derived in terms of the
$L^\infty$-norms of $\nabla u$ and $\nabla d$ in
\cite{Huang_Wang_Wen2011} and in terms of integral of
$L^\infty$-norms of $\nabla u$ and the BMO-norm of $\nabla d$,
in \cite{Chen_2012}. The global existence of weak solutions with large
initial data is still an outstanding open problem for  dimensions
larger then or equal to 3. So far, only results in one space
dimension have been obtained, for instance, we refer to
\cite{Ding2011, Ding2012}. In \cite{Hu_Wu_2012}, the existence and
uniqueness of global strong solutions to the Cauchy problem is
proved   for 3D in critical Besov spaces provided that the initial
data is close to an equilibrium state.

The overview of numerical methods used for the nematic liquid
crystal flows can be found in \cite{Badia_2011}.

At the same time,   up to now there are no results on the classical
solvability for the main initial and initial-boundary problems for
the full system (\ref{(1.1)}) -- (\ref{(1.5)}). In the present paper
we prove that a global classical solution to the Cauchy problem does
not exists, and that  this phenomenon is due to the presence of the
viscosity term. The nature of  loss of smoothness is similar to the
case of the compressible Navier-Stokes equation and the anisotropic
features do not influence  this phenomenon.

\section{A general blowup result}

In this paper we study  the system (\ref{(1.1)}) –- (\ref{(1.5)})
with the initial condition
\begin{equation}\label{ic}
(\rho, u, d)\Big|_{t=0} = (\rho_0, u_0, d_0),\quad d_0 \in {\mathcal
S}^2.
\end{equation}

We introduce the following natural functionals  defined on the
solution of the system (\ref{(1.1)}) -- (\ref{(1.5)}): mass
$$m(t)=\int\limits_{{\mathbb R}^3}\rho \, d x,$$  linear momentum
$$P(t)=\int\limits_{{\mathbb R}^3}\rho u\, d x,$$
and  total energy
\begin{equation}\label{energy}
\begin{aligned}
 E(t) =
 \frac12
\int\limits_{{\mathbb R}^n}\,\rho |u|^2\,dx +
\frac{1}{\gamma-1}\int\limits_{{\mathbb R}^3}\,p\,dx + \frac12
\int\limits_{{\mathbb R}^3}\,J\rho |\omega|^2\,dx  +
\int\limits_{{\mathbb R}^3}\,\rho W\,dx =\\
E_k(t)+E_i(t)+E_{p,1}(t) +E_{p,2}(t)\ge 0,
\end{aligned}
\end{equation}
where  $E_k(t)$, $E_i(t)$, and $E_{p,i},\,i=1,2,$ are the kinetic,
internal, and potential components of energy, respectively.

\begin{definition}
A solution $(\rho, u,d)$ to the Cauchy problem (\ref{(1.1)}) --
(\ref{(1.5)}), (\ref{ic}) belongs to the class $\mathfrak K$ if the
solution is classical, $\rho\ge 0$, the mass $m(t)$, linear momentum
$P(t)$, and total energy are finite for all $t\ge 0$, and, in addition, the
mass and linear momentum are conserved, i.e.,
$m(t)=m=\rm const, $
$P(t)=P=\rm const.$
\end{definition}

Thus, if the solution belongs to the class $\mathfrak K$, then
$$
\rho\in L^1\cap L^\gamma, \, \sqrt{\rho} u \in L^2,\, \sqrt{\rho}
\omega \in L^2,\, \nabla d\in L^2, \, (d,{\rm curl}\, d) \in L^1, \,
d \times {\rm curl}\, d \in L^1,
$$
where $L^r=L^r({\mathbb R}^3).$

Further, if we additionally assume
\begin{equation}\label{condition}
  u\in H^1,\quad d^T\,A\,d \in L^2, \quad Ad \in L^2, \quad N \in
L^2.
\end{equation}
A direct calculation with smooth solution $(\rho,\,u, \,d)$
to the system (\ref{(1.1)}) -- (\ref{(1.5)}) yields (cf.
\cite{Lin_Liu2000}, Theorem 2.1),
\begin{align}\label{deriv_energy}
 \frac{d}{dt}E(t) =& - \int\limits_{{\mathbb R}^3}\,
\left(\mu_1\,|d^T\,A\,d|^2 + \frac{\mu_4}{2} |\nabla u|^2 + (\mu_5 +
\mu_6) |Ad|^2\right)\,dx \nonumber \\
&+\lambda_1 \int\limits_{{\mathbb R}^3}\, |N|^2 \,dx  + (\lambda_2-
\mu_2 - \mu_3)\,\int\limits_{{\mathbb R}^3}\,(N,Ad)\,dx.
\end{align}
Here and throughout, we always assume that
\begin{equation}\label{lambda_mu}
\lambda_1 < 0,\quad \mu_5 + \mu_6 \ge 0, \quad \mu_1 \ge 0,\quad
\mu_4
> 0.
\end{equation}
 These
assumptions are necessary conditions for the dissipation of the
director field \cite{{Ericksen1991}, {Leslie1979}}.

Indeed, as was shown in \cite{full_system} by means of
the H\"older
inequality, if (\ref{lambda_mu}) holds, under the additional
assumption
\begin{equation}\label{additional}
|\lambda_2 -( \mu_2 + \mu_3)| \le 2\sqrt{-\lambda_1(\mu_5+\mu_6)}
\end{equation}
one has the following energy inequality: the basic energy law
(without Parodi's relation) holds
\begin{equation}\label{energy_law}
\frac{d} {dt} E(t) \le - \int\limits_{{\mathbb R}^3}\,
\left(\mu_1\,|d^T\,A\,d|^2 + \frac{\mu_4}{2} |\nabla
u|^2\right)\,dx \leq 0.
\end{equation}

Our main result is the following theorem.

\begin{theorem}\label{T}
Suppose that  $\gamma\ge \frac{6}{5}$ and inequalities
(\ref{lambda_mu}) and (\ref{additional}) hold. Then there is no
global in time solution to the Cauchy problem (\ref{(1.1)}) --
(\ref{(1.5)}), (\ref{ic}) in the class $\mathfrak K$ satisfying
(\ref{condition}).
\end{theorem}

The proof is similar to \cite{R_2008}, \cite{R_2009}, where a
blow-up result was proved for the compressible Navier-Stokes system
and the compressible magnetohydrodynamics system.

The next lemma is a key technical point of the  proof of the
theorem.

\begin{lemma} \label{lemma_stat}Let $\gamma\ge \frac{6}{5}$ and $u\in H^1({\mathbb R}^3)$. If $|P|\ne 0,$ then there
exists a positive constant $C$ such that for the solutions from the
class $\mathfrak K$ the following inequality holds:
\begin{equation}\label{lemma}
\int\limits_{{\mathbb R}^3}\,|Du|^2\, dx \ge
\,C\,E_i^{-\frac{1}{3(\gamma-1)}}(t).
\end{equation}
\end{lemma}

\begin{proof}
 First, from the inequality
\begin{equation}\label{embedding}
\left(\;\int\limits_{{\mathbb R}^n}\,
|u|^{\frac{2n}{n-2}}\,dx\right)^{\frac{n-2}{n}} \le C_1
\,\int\limits_{{\mathbb R}^n}\, |Du|^2 \,dx,
\end{equation}
where the constant $C_1>0 $ depends on $n,$  $\,n\ge 3, $
(\cite{Hebey}, p.22) we get $H^1({\mathbb R}^3)\subset L^6({\mathbb
R}^3)$. Then, using $\gamma \geq \frac{6}{5}$, the H\"older inequality gives
\begin{equation}\label{L1}
|P|=\left|\;\int\limits_{{\mathbb R}^3}\,\rho u \,dx\right|\le
\left(\;\int\limits_{{\mathbb R}^3}\,\rho^{\frac{6}{5}}
\,dx\right)^{\frac{5}{6}}
\left(\;\int\limits_{{\mathbb R}^3}\, |u|^6
\,dx\right)^{\frac{1}{6}}.
\end{equation}
Then, by the Jensen inequality (e.g., \cite{Kuczma_2009}, Theorem
8.1.3), we have for
 $\gamma\ge\frac{6}{5}$,
\begin{equation}\label{L2}
\left(\frac{1}{m}\,\int\limits_{{\mathbb R}^n}\,\rho^{\frac{6}{5}}
\,dx\right)^{5(\gamma-1)}\le \frac{\int\limits_{{\mathbb
R}^n}\,\rho^\gamma \,dx}{m}=\frac{(\gamma-1) E_i(t)}{mA}.
\end{equation}
Thus, (\ref{L1}) and (\ref{L2}) imply
\begin{equation}\label{L3}
|P|\le C_2\, \left(E_i(t)\right)^{\frac{1}{6(\gamma-1)}}
\left(\int\limits_{{\mathbb R}^3}\, |u|^6 \,dx\right)^{\frac{1}{6}},
\end{equation}
with the constant
$C_2=m^{\frac{5}{6}}\left(\frac{\gamma-1}{mA}\right)^{\frac{1}{6(\gamma-1)}}>0.$
Thus, the inequality (\ref{lemma}) follows from (\ref{embedding}),
(\ref{L3}),
 with the constant $C=\frac{|P|^2 C_2}{C_1^2}.$
\end{proof}

 \vskip1cm

{\bf Proof of Theorem \ref{T}}. Since $E_i(t)\leq E(t)\leq E(0)$ by \eqref{energy_law}, from
(\ref{additional}) and  (\ref{lemma})  we conclude
$$\dfrac {d}{dt}E(t)\le -\mu_4 \,C \,(E_i(t))^{-\frac{1}{3(\gamma-1)}}\le -\mu_4 \,C \,(E(0))^{-\frac{1}{3(\gamma-1)}}.  $$
which contradicts the non-negativity of $E(t)$ for all $t>0$. This concludes the proof of the theorem. $\Box$

\end{document}